\documentclass[12pt]{amsart}

\usepackage{amsmath,amssymb,epic,lscape}

\newtheorem{theorem}{Theorem}[section]

\theoremstyle{definition}

\theoremstyle{remark}

\numberwithin{equation}{section}
\providecommand{\bysame}{\leavevmode\hbox to3em{\hrulefill}\thinspace}

\def\DJ{{\hbox{D\kern-.8em\raise.15ex\hbox{--}\kern.35em}}}
\def\DJo{$\;$\kern-.4em
    \hbox{D\kern-.8em\raise.15ex\hbox{--}\kern.35em okovi\'c}}

\def\NSERC{Supported in part by an NSERC Discovery Grant.}

\def\bZ{{\mbox{\bf Z}}}

\def\pE{{\mathcal E}}

\def\Gr{{ G_{\rm BS} }}

\renewcommand{\subjclassname}{\textup{2000} Mathematics Subject
Classification }

\begin{document}

\title[Hadamard matrices from base sequences]
{Hadamard matrices from base sequences: An example}

\author[D.\v{Z}. \DJ okovi\'{c}]
{Dragomir \v{Z}. \DJ okovi\'{c}}

\address{Department of Pure Mathematics, University of Waterloo,
Waterloo, Ontario, N2L 3G1, Canada}

\email{djokovic@uwaterloo.ca}

\thanks{\NSERC}


\date{}

\begin{abstract}
There are several well-known methods that one can use to
construct Hadamard matrices from base sequences $BS(m,n)$. 
In view of the recent classification of base sequences
$BS(n+1,n)$ for $n\le30$, it may be of interest to show
on an example how prolific these methods are. For that
purpose we have selected the Hadamard matrices of order 60.
By using these methods and the transposition map we have 
constructed 1759 nonequivalent Hadamard matrices of order 60.
\end{abstract}

\maketitle
\subjclassname{ 05B20, 05B30 }
\vskip5mm

\section{Introduction} \label{Uvod}

Recall that a Hadamard matrix of order $m$ is a 
$\{\pm1\}$-matrix $H$ of size $m\times m$ such that $HH^T=mI_m$, 
where $T$ denotes the transpose and $I_m$ the identity matrix.
Let us denote by $H(m)$ the set of Hadamard matrices of order 
$m$. By abuse of language, we say that $H(m)$ exist if 
$H(m)\ne\emptyset$. If $m>2$ and $H(m)$ exist, then $m$ is 
divisible by 4. Two Hadamard matrices $A,B\in H(m)$ are said
to be {\em equivalent} if $B=PAQ$ for some signed permutation
matrices $P$ and $Q$.

By {\em binary} respectively {\em ternary sequence} we mean a 
sequence $A=a_1,a_2,\ldots,a_m$ whose terms belong to $\{\pm1\}$ 
respectively $\{0,\pm1\}$. 
To such a sequence we associate the polynomial 
$A(z)=a_1+a_2z+\cdots+a_mz^{m-1}$.
We refer to the Laurent polynomial $N(A)=A(z)A(z^{-1})$ as 
the {\em norm} of $A$.
{\em Base sequences} $(A;B;C;D)$ are quadruples of binary sequences,
with $A$ and $B$ of length $m$ and $C$ and $D$ of length $n$, 
and such that $N(A)+N(B)+N(C)+N(D)=2(m+n)$.
We denote the set of such base sequences by $BS(m,n)$.
We shall demonstrate that the base sequences and, their special 
cases, normal and near-normal sequences play an important role 
in the construction of Hadamard matrices \cite{HCD,SY}. 
The recent discovery of a Hadamard matrix of order 
428 \cite{KT1} used a $BS(71,36)$, constructed specially for 
that purpose.

As explained in \cite{DZ1}, we can view the normal sequences $NS(n)$ 
and near-normal sequences $NN(n)$ as subsets of $BS(n+1,n)$. 
For normal sequences $2n$ must be a sum of three squares, and 
for near-normal sequences $n$ must be even or 1.
The base sequences $(A;B;C;D)\in BS(n+1,n)$ are {\em normal}
respectively {\em near-normal} if $b_i=a_i$
respectively $b_i=(-1)^{i-1}a_i$ for all $i\le n$.

There are several known methods that one can use to
construct Hadamard matrices from base sequences $BS(m,n)$. 
In view of the recent classification \cite{DZ4} of base sequences
$BS(n+1,n)$ for $n\le30$, it may be of interest to show
on an example how prolific these methods are. For that
purpose we have selected the Hadamard matrices of order 60.
By using these methods we have constructed 1086 Hadamard matrices 
of order 60. Exactly 1012 of them are pairwise nonequivalent.
By taking transposes of these 1012 matrices, we obtain
additional 747 equivalence classes. Thus in total we have 
constructed 1759 equivalence classes of $H(60)$ by using
base sequences and the transposition. This is probably a very 
tiny portion of the totality of equivalence classes of $H(60)$. 
In that regard we mention that the incomplete classification, 
carried out very recently in \cite{KT2}, shows that the number 
of classes of $H(32)$ is well over 13 milions.

In Section \ref{BS87} we describe the construction of 
$H(8n+4)$ from $BS(n+1,n)$ and we summarize the results of our
computation in Table 1. The 558 pairwise nonequivalent Hadamard 
matrices constructed in this section are listed in Table 2.

Yang's paper \cite{CHY} contains four powerful ``multiplication
theorems''. The proofs of these theorems in {\it loc. cit.} 
are based on Yang's generalization of Lagrange's theorem on 
the sum of four squares to the ring of Laurent polynomials 
$\bZ[x,x^{-1}]$ with integer coefficients. 
This beautiful generalization deserves to be wider known. 
In a recent paper with K. Zhao \cite{DjZ} we have shown that 
Yang's generalization is essentially unique. 

As $15$ is a composite number, each of the four multiplication 
theorems can be used to construct base sequences $BS(15,15)$, 
and then construct Hadamard matrices of order 60. 
The results of these computations are described
in Sections 4-7. For the convenience of the reader, in each
of these four sections we state explicitly the multiplication 
theorem that we use. Three misprints in the statement of
Yang theorems have been observed in \cite{DZ3}.
The number of equivalence classes of $H(60)$ constructed in
these four sections are 192, 208, 64 and 64 respectively.
Their representative matrices are listed in 
the Appendix in Tables 3-6 respectively. 

In Section \ref{Kod} we describe the encoding of base sequences
$BS(n+1,n)$ that we use in Table 1, and have used in several
of our previous papers.

The appendix contains the Tables 2-6. We also explane there
how to interprete the entries of the tables.

\section{Preliminaries} \label{Oznake}

All computations were carried out in Magma \cite{Mag} modulo 
the tables of base sequences constructed in \cite{DZ1}.
In fact we only make use of Table 2 of that paper.
The main reason for using Magma is that it provides a (small) 
database of Hadamard matrices, in particular 256 matrices of 
order 60, and a useful collection of functions for working
with these matrices. The 1759 classes mentioned above are all different from these 256. The most valuable functions for us
were ``IsHadamard'' for testing whether a matrix is a
Hadamard matrix,  ``HadamardCanonicalForm'' which 
provides a test for equivalence of Hadamard matrices, and
``HadamardMatrixToInteger'' and its inverse 
``HadamardMatrixFromInteger'' for compact representation
of Hadamard matrices.

All Hadamard matrices in this note are constructed by using
the Goethals--Seidel array:
$$
\left[ \begin{array}{rrrr}
		Z_0 & Z_1R & Z_2R & Z_3R \\
		-Z_1R & Z_0 & -RZ_3 & RZ_2 \\
		-Z_2R & RZ_3 & Z_0 & -RZ_1 \\
		-Z_3R & -RZ_2 & RZ_1 & Z_0
\end{array} \right], $$
where $Z_0,Z_1,Z_2,Z_3$ are suitable circulant matrices,
and $R$ denotes the matrix having ones on the back-diagonal
and all other entries zero. 

Let us make a remark about this array. It is not hard to verify that
if we permute the circulants $Z_0,Z_1,Z_2,Z_3$ then
the new Hadamard matrix obtained from the above array will
be equivalent to the original one provided that the 
permutation is even. If it is odd then the two Hadamard 
matrices may be nonequivalent. 
This is used in Section \ref{Yang3}.

We now list additional notation and definitions that we need.

We separate the sequences by a semicolon, and
use the comma as the concatenation operator.
The symbol $0_s$ denotes the sequence of $s$ zeros.
For a sequence $A=a_1,a_2,\ldots,a_m$ we denote by $A'$ 
the reversed sequence, i.e., $A'=a_m,a_{m-1},\ldots,a_1$. 
Thus we have $a'_k=a_{m+1-k}$ for $k=1,2,\ldots,m$. 
If $f\in\{+1,-1\}$ then $fA$ is the sequence
$fa_1,fa_2,\ldots,fa_m$.
For sequences $A$ and $B$ of length $n$, $A\pm B$ denotes the 
sequence with terms $a_i\pm b_i$, $i=1,2,\ldots,n$.
For sequences $A=a_1,a_2,\ldots,a_{m+1}$ and $C=c_1,c_2,\ldots,c_m$
we denote by $A/C$ the interlaced sequence
$$
A/C=a_1,c_1,a_2,c_2,\ldots,a_m,c_m,a_{m+1}.
$$

We say that two ternary sequences $G$ and $H$ of length $n$ are 
{\em disjoint} if at most one of $g_i$ and $h_i$ is nonzero
for each index $i$. We recall that $T$-{\em sequences} are 
quadruples $(A;B;C;D)$ of pairwise disjoint ternary sequences of 
length $n$ such that $N(A)+N(B)+N(C)+N(D)=n$.  
We denote by $TS(n)$ the set of $T$-sequences of length $n$. 
It has been conjectured that $TS(n)\ne\emptyset$ for all 
integers $n\ge1$, and it is known that this is true for 
$n\le100$ different from $79$ and $97$. There is a map
\begin{equation} \label{Kons-3}
TS(n) \to BS(n,n)
\end{equation}
sending $(A;B;C;D) \to (Q,R,S,T)$ where
\begin{eqnarray*}
	Q &=& A+B+C+D;\\ R &=& A+B-C-D;\\
	S &=& A-B+C-D;\\ T &=& A-B-C+D.
\end{eqnarray*}

\section{Construction of $H(60)$ from $BS(8,7)$} 
\label{BS87}

For any finite binary sequence $X$ let $Z_X$ denote the
circulant matrix having $X$ as its first row.
It is well known (see e.g. \cite{SY}) that there is a map
\begin{equation} \label{Kons-1}
BS(d,d) \to H(4d).
\end{equation}
sending $(M;U;V;W)$ to the Hadamard matrix $H$ constructed by plugging in the circulants $Z_M,Z_U,Z_V,Z_W$ for $Z_0,Z_1,Z_2,Z_3$ 
in the Goethals--Seidel array.

One can also use base sequences $BS(m,n)$ with $m$ and $n$
arbitrary to construct Hadamard matrices of order $4(m+n)$.
For that purpose we just compose the above map with the map
$$ BS(m,n) \to BS(m+n,m+n) $$
which sends $(A;B;C;D)\to (A,C;A,-C;B,D;B,-D)$.
In particular, for $m=n+1$, we obtain the map
\begin{equation} \label{Kons-2}
BS(n+1,n) \to H(8n+4).
\end{equation}

In our recent paper \cite{DZ4} we have introduced a new
equivalence relation in $BS(n+1,n)$ which the reader should
consult for further details. As this relation is
not easy to describe, we shall just say that there is
a group $\Gr$ of order $2^{12}$, whose definition
depends on the parity of $n$, which acts on $BS(n+1,n)$
so that its orbits are exactly the equivalence classes.
By using this group, it is easy to generate in Magma
the whole equivalence class from its representative.
We point out that the map (\ref{Kons-2}) may produce many
nonequivalent Hadamard matrices from a single equivalence
class of base sequences (see Table 1).

Let $\pE\subseteq BS(n+1,n)$ be an equivalence class. From
above it follows that the cardinality of $\pE$ must
be a power of two, $2^k$ with $k\le12$. We are interested in
the case $n=7$ in which case $2n+1=15$ and $8n+4=60$. The
set $BS(8,7)$ splits into 17 equivalence classes. 

The results of our computation in this case are summarized
in Table 1. The listing of the equivalence classes of
$BS(8,7)$ is taken from \cite[Table 2]{DZ1}. The second
column of Table 1 lists the representatives $(A;B;C;D)$ of 
the equivalence classes $\pE$ of $BS(8,7)$. These 
representatives are written in encoded form (see  
Section \ref{Kod} for the definition and our conventions for this 
encoding). For each representative we record the cardinality 
$\#\pE$ of $\pE$ and the number of equivalence classes 
of $H(60)$ constructed from $\pE$ via (\ref{Kons-2}). 

It turns out that any two nonequivalent base sequences in
$BS(8,7)$ produce two nonequivalent Hadamard matrices.
We do not know whether this is true in general.

{\bf Question :} Is it true that two nonequivalent base
sequences in $BS(n+1,n)$ always produce via (\ref{Kons-2}) 
two nonequivalent Hadamard matrices of order $8n+4$?

The sum of the numbers in the last column of Table 1 is 558. 
Consequently, we have constructed 558 equivalence
classes of $H(60)$. 

\begin{center}
\begin{tabular}{rrrc}
\multicolumn{4}{c}{Table 1: $H(60)$ from $BS(8,7)$} \\ \hline 
\multicolumn{1}{c}{} & \multicolumn{1}{c}{$AB;\,CD$} & 
\multicolumn{1}{c}{$\#\pE$} &
\multicolumn{1}{c}{\#Had.} \\ \hline
1 & 0165;\, 6123 & $2048$ & $32$ \\ 
2 & 0165;\, 6141 & $4096$ & $64$ \\ 
3 & 0166;\, 6122 & $2048$ & $32$ \\ 
4 & 0173;\, 6161 & $2048$ & $32$ \\ 
5 & 0173;\, 6411 & $4096$ & $64$ \\ 
6 & 0183;\, 6121 & $2048$ & $32$ \\ 
7 & 0613;\, 1623 & $2048$ & $32$ \\ 
8 & 0614;\, 1641 & $4096$ & $64$ \\ 
9 & 0615;\, 1263 & $2048$ & $32$ \\ 
10 & 0615;\, 1272 & $512$ & $8$ \\ 
11 & 0616;\, 1262 & $2048$ & $32$ \\ 
12 & 0618;\, 1261 & $2048$ & $32$ \\ 
13 & 0635;\, 1621 & $1024$ & $16$ \\ 
14 & 0638;\, 1620 & $2048$ & $32$ \\ 
15 & 0641;\, 1622 & $2048$ & $32$ \\ 
16 & 0646;\, 1222 & $256$ & $6$ \\ 
17 & 0646;\, 1260 & $1024$ & $16$ \\ 
\hline
\end{tabular} \\
\end{center}

\section{H(60) from Yang's Theorem 1} \label{Yang1}

The theorem \cite[Theorem 1]{CHY} gives a map
\begin{equation} \label{Yang-1}
NS(n) \times BS(s,t) \to TS(d),\quad d=(2n+1)(s+t).
\end{equation}
Normal sequences $NS(n)$ can be written as 
$(F,+;F,-;G+H;G-H)$ or $(F,-;F,+;G+H;G-H)$, 
where $F$, $G$ and $H$ are uniquely determined sequences 
of length $n$ such that $F$ is binary while $G$ and $H$
are ternary and disjoint. The map (\ref{Yang-1})
is defined in terms of the sequences $F$, $G$ and $H$.
This map sends the ordered pair, having these
normal sequences as the first component and 
$(A;B;C;D)\in BS(s,t)$ as the second, to the $T$-sequences
$(Q,R,S,T)$ where
\begin{equation} \label{QRST}
\left[ \begin{array}{c} Q\\R\\S\\T \end{array} \right] =
\left[ X_1,X_2, \ldots ,X_n,X_{n+1} \right]
\end{equation}
and the blocks $X_k$, $k=1,2,\ldots,n$ and $X_{n+1}$ 
are given by
$$
X_k=\left[ \begin{array}{llr} 
f'_kA, & g_kC+h_kD, & 0_{s+t} \\
f'_kB, & -h'_kC+g'_kD, & 0_{s+t} \\
0_{s+t}, & g'_kA-h_kB, & -f_kC \\
0_{s+t}, & h'_kA+g_kB, & -f_kD \end{array}\right], \quad
X_{n+1}=\left[ \begin{array}{rr} 
-B', & 0_t \\
A', & 0_t \\
0_s, & -D' \\
0_s, & C' \end{array} \right].
$$

Two misprints in the expression for $\tau_k$ in \cite[p. 770]{CHY} 
have been corrected. Instead of our sequences $g'_kA-h_kB$ and 
$h'_kA+g_kB$, inside the block $X_k$, Yang has $g'_kA-h'_kB$ and 
$h'_kA+g'_kB$, respectively. At a first glance our sequences appear 
to be in error since these have to be ternary. In fact they
are binary since Theorem \ref{KKS-teorema}
guarantees that $g_k=0$ iff $g'_k=0$ and $h_k=0$ iff $h'_k=0$.
Hence, exactly one of $g'_k$ and $h_k$ is zero and the same
is true for $h'_k$ and $g_k$.

There exists only one equivalence class of base sequences
$BS(2,1)$ and the same is true for $BS(3,2)$. Their representatives
are $03;1$ and $0;0$ in encoded form, or explicitly
\begin{eqnarray*}
0;0 &=&  +,+;\, +,-;\, +;\, + \\
03;1 &=& +,-,+;\, +,-,-;\, +,+;\, +,+.
\end{eqnarray*}
These are also representatives of the equivalence classes
of normal sequences $NS(1)$ and $NS(2)$ respectively.

By applying the above theorem, with $s=3$, $t=2$ and $n=1$, we
compute the image of the whole set $NS(1)\times BS(3,2)$ and 
then apply the map (\ref{Kons-2}) to this image.
We obtain 128 equivalence classes of $H(60)$. Another 64
equivalence classes are obtained by taking $s=2$, $t=1$ and $n=2$,
i.e., bu using the set $NS(2)\times BS(2,1)$.
These $128+64=192$ equivalence classes turn out to be distinct.

\section{H(60) from Yang's Theorem 2} \label{Yang2}

The theorem \cite[Theorem 2]{CHY} gives a map
\begin{equation} \label{Yang-2}
NS(n) \times BS(s,t) \to BS(d,d),\quad d=n(s+t).
\end{equation}
We write normal sequences $NS(n)$ as in the previous section,
i.e., as $(F,+;F,-;G+H;G-H)$ or $(F,-;F,+;G+H;G-H)$.
The map (\ref{Yang-2}) sends the ordered pair consisting of these 
normal sequences  and the base sequences $(A;B;C;D)\in BS(s,t)$ 
to the base sequences $(Q,R,S,T)$ defined by
$$
\left[ \begin{array}{c} Q\\R\\S\\T \end{array} \right] =
\left[ X_1,X_2, \ldots ,X_n \right], \quad
X_k=\left[ \begin{array}{lr} 
f'_kA, & g_kC+h_kD \\
f'_kB, & -h'_kC+g'_kD \\
g'_kA-h_kB, & -f_kC \\
h'_kA+g_kB, & -f_kD \end{array} \right].
$$

There are two possibilities. First we take $s=3$, $t=2$ and $n=3$.
There exists only one equivalence class of normal sequences
$NS(3)$. As its representative we can take
$$ 06;11=+,+,-,+;\, +,+,-,-;\, +,+,+;\, +,-,+. $$
By computing the image of $NS(3) \times BS(3,2)$ under the 
map (\ref{Yang-2}) and applying the map (\ref{Kons-2}), 
we obtain 80 equivalence classes of $H(60)$. 

The second possibility is to take $s=2$, $t=1$ and $n=5$.
Again there exists only one equivalence class of normal sequences
$NS(5)$. As its representative we can take $016;640$, i.e.,
$$ +,+,+,-,+,+;\, +,+,+,-,+,-;\, +,+,+,-,-;\, +,-,+,+,-. $$
In this case we obtain 128 equivalence classes of $H(60)$.

These $80+128=208$ equivalence classes turn out to be distinct.

\section{H(60) from Yang's Theorem 3} \label{Yang3}

The theorem \cite[Theorem 3]{CHY} gives a map
\begin{equation} \label{Yang-3}
NN(n) \times BS(s,t) \to TS(d),\quad d=(2n+1)(s+t),
\end{equation}
where $n=2m$ is even.
To describe this map, we shall write near-normal sequences in
$NN(n)$ in the form
$$ ((Y,+)/X;(Y,-)/(-X);G+H;G-H), $$
where $X$ and $Y$ are binary sequences and $G$ and $H$ are disjoint 
ternary sequences, all of length $n$. 
This map sends the ordered pair, having these
near-normal sequences as the first component and 
$(A;B;C;D)\in BS(s,t)$ as the second, to the $T$-sequences
$(Q;R;S;T)$ where
\begin{eqnarray*}
\left[ \begin{array}{c} Q\\R\\S\\T \end{array} \right] &=&
\left[ \begin{array}{c}
U_1,U_2, \ldots ,U_m,U_{m+1} \\
V_{m+1},V_m, \ldots ,V_2,V_1 \end{array} \right],
\end{eqnarray*}
the blocks $U_k$ and $V_k$, $k\le m$, are given by
\begin{eqnarray*}
U_k &=& \left[ \begin{array}{lllr} 
g'_{2k-1}A-h_{2k-1}B, &-y_kC, &g'_{2k}A-h_{2k}B, &-x_kD' \\
h'_{2k-1}A-g_{2k-1}B, &-y_kD, &h'_{2k}A+g_{2k}B, &x_kC'
\end{array} \right], \\
V_k &=& \left[ \begin{array}{rlll} 
-x_kB, &g_{2k}C'+h_{2k}D', &y'_kA', &g_{2k-1}C'+h_{2k-1}D' \\
x_kA, &g'_{2k}D'-h'_{2k}C', &y'_kB', &g'_{2k-1}D'-h'_{2k-1}C'
\end{array} \right],
\end{eqnarray*}
and
$$
U_{m+1}= \left[ \begin{array}{lrl} 
0_s, & -D', & 0_{n(s+t)} \\
0_s, & C', & 0_{n(s+t)} \end{array} \right], \quad
V_{m+1}= \left[ \begin{array}{lrl} 
0_{n(s+t)}, & -B, & 0_t \\
0_{n(s+t)}, & A, & 0_t \end{array} \right].
$$

There exists only one equivalence class in $NN(2)$. As its 
representative we can take
$$ 02;1=+,-,+;\, +,+,-;\, +,+;\, +,+. $$
After computing the image of the map (\ref{Yang-3}) with 
$s=2$, $t=1$ and $n=2$ and applying the maps (\ref{Kons-3}) and
(\ref{Kons-2}) in succession, we obtain 32 equivalence classes 
of $H(60)$. Another 32 equivalence classes are obtained by 
swaping the first two components of the base sequences produced 
by the map (\ref{Kons-3}) (see the remark made in 
Section \ref{Oznake}).
These $32+32=64$ equivalence classes turn out to be distinct.

\section{H(60) from Yang's Theorem 4} \label{Yang4}

The theorem \cite[Theorem 4]{CHY} gives a map
\begin{equation} \label{Yang-4}
BS(m+1,m) \times BS(n+1,n) \to BS(d,d),\quad d=(2m+1)(2n+1),
\end{equation}
which sends the ordered pair $\left( (A;B;C;D),(F;G;H;E) \right)$ 
to $(Q;R;S;T)$ defined again by the formula (\ref{QRST}) but 
the blocks $X_k$, $k\le n$, and $X_{n+1}$ are now given by
$$
X_k= \left[ \begin{array}{ll} 
f'_kA/g_kC, & (-e_kB')/h_kD \\
f'_kB/g'_kD, & e_kA'/(-h_kC) \\
g'_kA/(-f_kC), & (-h_kB)/(-e_kD') \\
g_kB/(-f_kD), & h'_kA/e_kC' \end{array} \right], \quad
X_{n+1}= \left[ \begin{array}{l} 
f_1A/g'_1C \\
f_1B/g_1D \\
g_1A/(-f'_1C) \\
g'_1B/(-f'_1D) \end{array} \right]
$$
(a misprint in the expression for $\beta_k$
in \cite[p. 773]{CHY} has been corrected).

We apply this theorem, with $m=1$ and $n=2$. We compute the
image of $BS(2,1)\times BS(3,2)$ and then apply the map 
(\ref{Kons-2}). We obtain 32 equivalence classes of $H(60)$. 
Another 32 equivalence classes are obtained by using the set
$BS(2,1) \times BS(3,2)$. These $32+32=64$ equivalence classes turn
out to be distinct.

\section{The encoding scheme} \label{Kod}

Let $(A;B;C;D) \in BS(n+1,n)$. For $n$ even (odd) set 
$n=2m$ ($n=2m+1$). Decompose $(A;B)$ into quads
$$ \left[ \begin{array}{ll} a_i & a_{n+2-i} \\ 
b_i & b_{n+2-i} \end{array} \right],\quad i=1,2,\ldots,
\left[ \frac{n+1}{2} \right], $$ 
and, if $n$ is even, the central column
$ \left[ \begin{array}{l} a_{m+1} \\ b_{m+1} \end{array} \right]. $
Similar decomposition is valid for $(C;D)$.
The quad encoding is based on \cite[Theorem 1]{KKS}.
\begin{theorem} \label{KKS-teorema}
The sum of the four quad entries is $2 \pmod{4}$ for the
first quad of $(A;B)$ and is $0 \pmod{4}$ for all other 
quads of $(A;B)$ and for all quads of $(C;D)$.
\end{theorem}

There are 8 possibilities for the first quad of $(A;B)$:
\begin{eqnarray*}
1'=\left[ \begin{array}{ll} - & + \\ + & + \end{array} \right],\quad 
2'=\left[ \begin{array}{ll} + & - \\ + & + \end{array} \right],\quad 
3'=\left[ \begin{array}{ll} + & + \\ + & - \end{array} \right],\quad 
4'=\left[ \begin{array}{ll} + & + \\ - & + \end{array} \right], \\
5'=\left[ \begin{array}{ll} + & - \\ - & - \end{array} \right],\quad 
6'=\left[ \begin{array}{ll} - & + \\ - & - \end{array} \right],\quad 
7'=\left[ \begin{array}{ll} - & - \\ - & + \end{array} \right],\quad 
8'=\left[ \begin{array}{ll} - & - \\ + & - \end{array} \right].
\end{eqnarray*}
The possibilities for the remaining quads 
of $(A;B)$ and the quads of $(C;D)$ are:
\begin{eqnarray*}
1=\left[ \begin{array}{ll} + & + \\ + & + \end{array} \right],\quad 
2=\left[ \begin{array}{ll} + & + \\ - & - \end{array} \right],\quad 
3=\left[ \begin{array}{ll} - & + \\ - & + \end{array} \right],\quad 
4=\left[ \begin{array}{ll} + & - \\ - & + \end{array} \right], \\
5=\left[ \begin{array}{ll} - & + \\ + & - \end{array} \right],\quad 
6=\left[ \begin{array}{ll} + & - \\ + & - \end{array} \right],\quad 
7=\left[ \begin{array}{ll} - & - \\ + & + \end{array} \right],\quad 
8=\left[ \begin{array}{ll} - & - \\ - & - \end{array} \right].
\end{eqnarray*}
Finally, there are 4 possibilities for the central column:
$$
0=\left[ \begin{array}{l} + \\ + \end{array} \right],\quad
1=\left[ \begin{array}{l} + \\ - \end{array} \right],\quad
2=\left[ \begin{array}{l} - \\ + \end{array} \right],\quad
3=\left[ \begin{array}{l} - \\ - \end{array} \right].
$$

We encode $(A;B)$ by the symbol sequence
$p_1p_2 \ldots p_m p_{m+1}$,
where $p_i$ is the label of the $i$th quad except that for
$n$ even $p_{m+1}$ is the label of the central column. 
Similarly, we encode $(C;D)$ by $ q_1q_2 \ldots q_m$ respectively
$q_1q_2 \ldots q_m q_{m+1}$ when $n$ is even respectively odd. 
Here $q_i$, $i\le m$, is the label of the $i$th quad and, for
$n$ odd, $q_{m+1}$ is the label of the central column.

As an example, the base sequences
\begin{eqnarray*}
A &=& +,+,+,+,-,-,+,-,+; \\
B &=& +,+,+,-,+,+,+,-,-; \\
C &=& +,+,-,-,+,-,-,+; \\
D &=& +,+,+,+,-,+,-,+;
\end{eqnarray*}
are encoded as $3'6142;\, 1675$. In Table 1 and elsewhere in the
text we write 0 instead of $3'$.

\section{Appendix} \label{Dodatak}

The following tables contain the list of Hadamard matrices
constructed by the methods explained in Sections 3-7.
Since these matrices are constructed by using the Goethals--Seidel
array, i.e., via the map (\ref{Kons-2}), they can be stored very efficiently. Indeed it suffices to list only the base sequences 
$(A;B;C;D) \in BS(15,15)$ used to construct the Hadamard matrix.
We first concatenate these four constituent binary sequences to
obtain the binary sequence $A,B,C,D$ of length 60. Next we
replace each $-1$ term in this sequence with 0 to obtain a
$\{0,1\}$-sequence, say $S$, of length 60. Next we split $S$
into 15 pieces of length 4 each. Finally, we interprete each
piece as the binary representation of a hexadecimal digit
$0,1,\ldots,9,a,b,c,d,e,f$. Thus we obtain a sequence, $X$, of 
exactly 15 hexadecimal digits. Each entry in our tables is such 
a sequence. Clearly, one can easily reconstruct the base sequences
$(A;B;C;D)$ from $X$. 

Here is an example. We take $X=0dc41a77adbf5c8$, the first
hexadecimal sequence in Table 2. Each hexadecimal digit has
to be replaced by its binary representation by using exactly
four binary digits. Thus the hexadecimal digit 0 gets replaced
by the sequence $0,0,0,0$. Next the hexadecimal digit ``d''
is replaced by the sequence $1,1,0,1$ etc. We thus
obtain the sequence $S$ as
$$000011011100010000011010011101111010110110111111010111001000.$$
Finally, we replace each 0 with $-1$ and read off the four
binary sequences:
\begin{eqnarray*}
A &=& -,-,-,-,+,+,-,+,+,+,-,-,-,+,-; \\
B &=& -,-,-,-,+,+,-,+,-,-,+,+,+,-,+; \\
C &=& +,+,+,-,+,-,+,+,-,+,+,-,+,+,+; \\
D &=& +,+,+,-,+,-,+,+,+,-,-,+,-,-,-. \\
\end{eqnarray*}
To obtain the Hadamard matrix $H\in H(60)$, it remains to form
the circulants $Z_A$, $Z_B$, $Z_C$, $Z_D$ and plug them
(in that order) into the Goethals--Seidel array.

The Hadamard matrices within each table are pairwise nonequivalent.
The same is true for the list of all Hadamard matrices made up
from Tables 3, 5 and 6. However there is an overlap (32 matrices
in total) between these tables and Table 4. Namely, 24 matrices 
in Table 4 are equivalent to some matrices in Table 3 and another 
8 matrices in Table 4 are equivalent to some matrices in Table 5. 
Similarly, there is an overlap (42 matrices in total) 
between Table 2 and the other four tables. Thus these five tables
together contain representatives of 1012 equivalence classes
of $H(60)$.

\newpage

\begin{center}
\begin{tabular}{llll}
\multicolumn{4}{c}{Table 2: $H(60)$ from $BS(8,7)$} \\  \hline 
0dc41a77adbf5c8& a73b4f89f643eb7& f2ede4275f16b9d& b03b618bac4f5f6\\
e547cb71f44bef6& d7edae2434e86e2& d7c5ae753e427b7& ebc5d676c245837\\ 
411283db96e72a3& a7474f7105ba0c8& be6f7d23971729d& 7db8fa8f97b7289\\
eb3bd7893c4a7f6& 82ed04269ced3e2& b03b618bae475b7& 0d6e1b235ceebe2\\ 
ebedd6253d1a7dc& e56fcb2106e20a3& 826f0523971729d& b09160df5ee6ba3\\
d7c5ae77ca47937& 8247057396472b7& 413a838a9dbd3c8& ebedd6253f1279d\\ 
f23be58bafb7589& d7edae2436e06a3& e5c5ca75f5bbec8& 7db8fa8e9fb5389\\
e5edca2506e20a3& b0ed60275d1ebdc& 7d3afb8b97b7289& e591cadd06e20a3\\ 
a7894eedf793e8d& 4176831396c72a7& a7994ecdf7d3e85& eb1bd7cac08d86e\\
a7cd4e66f82df7a& eb0bd7eac0cd866& 4166833396872af& a7dd4e46f86df72\\ 
be4f7d6395df2c4& be5f7d43959f2cc& e54fcb610772091& e55fcb410732099\\
f233e59b5fd6b85& f223e5bb5f96b8d& e5e5ca35f773e91& 4f229fbb5f96b8d\\ 
e5f5ca15f733e99& 4f329f9b5fd6b85& 82f504169f35399& 82e504369f75391\\
4f669f335e86baf& 4f769f135ec6ba7& beb17c9e9e253bb& bea17cbe9e653b3\\ 
82a104be9f35399& 82b1049e9f75391& e55fcb4104ca0e6& d723afbac13d858\\
d733af9ac17d850& e54fcb61048a0ee& d71bafc93c8a7ee& d70bafe93cca7e6\\ 
824f05639777291& f24fe5635ddebc4& 4f989ecf5e86baf& 4f889eef5ec6ba7\\
825f05439737299& f25fe5435d9ebcc& e5f5ca15f4cbee6& e5e5ca35f48beee\\ 
b05f61435d9ebcc& b04f61635ddebc4& e50bcbe9f663eb3& e51bcbc9f623ebb\\
a75f4f42099c14c& a74f4f6209dc144& eb33d79837d0685& eb23d7b8379068d\\ 
eba1d6bfc8cf966& ebb1d69fc88f96e& 0d4e1b63addf5c4& 0da01abf5ccebe6\\
bef57c169d9d3cc& bee57c369ddd3c4& d777af13cac7927& f25fe5435e66bb3\\ 
eb67d73037d0685& eb77d710379068d& be777d1396c72a7& d7b1ae9c35d86c4\\
be677d3396872af& d7a1aebc35986cc& d7ddae453f9278d& d7cdae653fd2785\\ 
5946b3729ff5381& 4f569f533fb6789& a77f4f0194eb2e2& f313e7dac15d854\\
59c4b27796072bf& cf6f9f2035586d4& e5b9ca8e6a04d3f& cfed9e253f52795\\ 
a7034ff9651acdc& a7474f729bf5301& cfa99eac36406b7& 597eb30396e72a3\\
cf579f53cbb7909& b0d56057cdbf9c8& a6b94c8f940f2fe& 4f569f53cc4f9f6\\ 
e57fcb0196e32a3& 0d6e1b233cae7ea& b02b61ab3dbe7c8& 4fec9e273cae7ea\\
f2ede4273ea67ab& cfd59e56c1bd848& cfa99eac35b86c8& a7b94e8d65facc0\\ 
f391e6dc34a86ea& b09160df3f56795& a6814cff971729d& e53bcb8a9bf5301\\
a77f4f0298ed362& 65fcca069ced3e2& e547cb726a04d3f& a63b4d8a9e053bf\\ 
d3b1a69e202c47a& d3a1a6be206c472& bbe57636782cf7a& bbf57616786cf72\\
bbb1769e197c350& 61b0c29eee85daf& bba176be193c358& 61a0c2beeec5da7\\ 
dda1babd879308d& ddb1ba9d87d3085& 9e4f3d637686eaf& 9e5f3d4376c6ea7\\
dde5ba3586830af& ddf5ba1586c30a7& 11a022bf4f9698d& 11b0229f4fd6985\\ 
88a110bf2c6e5f2& 88b1109f2c2e5fa& 88a110bf4d3e9d8& 88b1109f4d7e9d0\\
cba196bc44688f2& cbb1969c44288fa& cbb1969c47d0885& cba196bc479088d\\ 
d30ba7e9df93b8d& d31ba7c9dfd3b85& ee4fdd632fd6585& ee5fdd432f9658d\\ 
\hline
\end{tabular} \\
\end{center}

\begin{center}
\begin{tabular}{llll}
\multicolumn{4}{c}{Table 2 (continued)} \\  \hline 
dda1babd853b0d8& ddb1ba9d857b0d0& ddf5ba15e6c3ca7& dde5ba35e683caf\\ 
dd17bbd2190c35e& dd07bbf2194c356& cb0797f045a88ca& cbc1967fbab7729\\
cb1797d045e88c2& 887b110b4c5e9f4& cbd1965fbaf7721& cb8596f620bc468\\ 
886b112b4c1e9fc& cb9596d620fc460& dd3fbb821ab4329& dd2fbba21af4321\\
bb6b772a1b0431f& bb7b770a1b44317& bbe9762e79ecf42& 9e533d5b74feee0\\ 
bbf9760e79acf4a& dd7bbb0a7b44f17& 9e433d7b74beee8& bb3f7782794cf56\\
dd6bbb2a7b04f1f& bb7b770a7b44f17& bb2f77a2790cf5e& bb6b7729e703c9f\\ 
bb7b7709e743c97& 881711d34dee9c2& bb6b772a7b04f1f& 880711f34dae9ca\\
88bd10872cbe5e8& 88ad10a72cfe5e0& 77bcee872cbe5e8& 77aceea72cfe5e0\\ 
776aef2b2c1e5fc& 777aef0b2c5e5f4& 77e8ee2f4ef69a1& 77f8ee0f4eb69a9\\
eed1dc5f4dee9c2& 116a232b4cfe9e0& cbf9960fbab7729& eec1dc7f4dae9ca\\ 
117a230b4cbe9e8& cbe9962fbaf7721& 6116c3d2eef5da1& 6106c3f2eeb5da9\\
61e8c22f75eeec2& 111623d32d0e5de& 9ead3ca774feee0& 61f8c20f75aeeca\\ 
110623f32d4e5d6& dd07bbf2794cf56& 9ebd3c8774beee8& dd17bbd2790cf5e\\
dd07bbf21ab4329& cb7b970847a088b& dd17bbd21af4321& cb6b972847e0883\\ 
d385a6f7b85f774& d395a6d7b81f77c& dd6bbb2a78fcf60& dd7bbb0a78bcf68\\
9e2f3da2ededdc2& 9e3f3d82edaddca& ee3fdd834d4e9d6& ee2fdda34d0e9de\\ 
61d4c257e647cb7& cafd94074cee9e2& 9e133dda7f54f95& acc558774e069bf\\
d347a772b3f5601& ca8194ff4f1699d& 53b8a68f4dfe9c0& cb479770d601abf\\ 
ca4795734ff6981& 9fd53e55e5bbcc8& 9f133fda7b54f15& 53c4a6774c0e9fe\\
9f573f51e643cb7& 537ea7032cee5e2& cb0397fab0ed662& d3b9a68d4e029bf\\ 
9fd53e5585bb0c8& f92bf3a987b3089& cb7f9700d711a9d& 9f2b3faa1bb4309\\
9e2b3daa7dbcfc8& 792af3abe5bfcc8& 86130ddbe757c95& cbfd9606b2e5623\\ 
ac4759734c0e9fe& 9fed3e2586a30ab& ac0359fb4d1e9dc& cbb9968f280e57e\\
9f573f527a44f37& 792af3aa7e44fb7& 79d4f257e44fcf6& 612ac3abe7b7c89\\ 
ed87daf0dd41bd7& b7eb6e28dcf1be1& 133c2786dfadb8a& ed15dbd7b1e7643\\
47688f2ee415cfd& 9dd33a5877e8e82& 9d053bf6e0b5c69& ec51d95fb4f76e1\\ 
b8bf708274b4ee9& b99772d2e3edc02& b92d73a78aff120& ed79db0c4eb89a8\\
9dfb3a091da23cb& b92d73a6e105c5f& 1dd23a5a76fcea0& c8eb902adf0db9e\\ 
b7796f0f2146457& 131427d7b61f6bc& edafdaa3b21f63c& b7876ef04eb89a8\\
b7516f5cde19bbc& e2fbc40a765ceb4& b7516f5cdf09b9e& edafdaa3b0f7661\\ 
b9fb720ae0b5c69& b941737f89a714b& b9bf72838b4f116& 472c8fa676fcea0\\
b7516f5c4f0899e& edc3da78deb9ba8& 47048ff6e4b5ce9& b8fb700ae5a5ccb\\ 
13b62692bc757f1& 13a626b2bc357f9& b9f37218172829a& b9e372381768292\\
af495f6fb077671& af595f4fb037679& 9d8f3ae2836d012& a059414ee585ccf\\ 
a049416ee5c5cc7& f5cbea684e789b0& 9d9f3ac2832d01a& 1360273ebcd57e5\\
f5dbea484e389b8& 1370271ebc957ed& a0f3401ae69dcac& a0e3403ae6ddca4\\ 
05a60ab2e435cf9& 05b60a92e475cf1& faf3f41a7564ed3& fae3f43a7524edb\\
c8b79093d63fab8& c8a790b3d67fab0& a0db404a7474ef1& a0cb406a7434ef9\\ 
\hline
\end{tabular} \\
\end{center}

\begin{center}
\begin{tabular}{llll}
\multicolumn{4}{c}{Table 2 (continued)} \\  \hline 
37246fb6be3d7b8& 37346f96be7d7b0& 9de33a3a836d012& 9df33a1a832d01a\\
9df33a1814d02e5& 9de33a3814902ed& 05e20a3a776ce92& 05f20a1a772ce9a\\ 
9db73a93e9c7d47& 9da73ab3e987d4f& 1348276ebc757f1& 1358274ebc357f9\\
9d9f3ac014d02e5& c8b79093d78fa8e& c8a790b3d7cfa86& 9d8f3ae014902ed\\ 
afcb5e684fc8986& 9d353b9417c8286& ec8fd8e3d6dfaa4& 9d253bb4178828e\\
afdb5e484f8898e& ec9fd8c3d69faac& 9d353b9416782b0& 9d253bb416382b8\\ 
9d593b4c16782b0& c8a790b3d437af9& 9d493b6c16382b8& c8b79093d477af1\\
1370271ebedd7a4& 1360273ebe9d7ac& af1d5fc44f68992& af0d5fe44f2899a\\ 
eca7d8b3d7cfa86& ecb7d893d78fa8e& 378e6ee2bd257db& 379e6ec2bd657d3\\
9d593b4c17c8286& 9d493b6c178828e& f5b7ea904c709f1& f5a7eab04c309f9\\ 
b78b6eeb433e818& edcfda602d705d1& 1d203bbe85950cd& e2f7c41285350d9\\
ed8bdaeb4066873& 85f70a12e135c59& ed19dbcc2d705d1& 0b8a16eadc65bf3\\ 
2f5c5f46dcc5be7& e24dc56684850ef& e209c5ee15942cd& b7e76e302e285ba\\
d075a116df3db98& f419e9ceddd5bc5& eda3dab82cc05e7& 854d0b66e3ddc04\\ 
85090bed1e6a3b2& ed8bdae8bf39798& a1b342991fda384& 474c8f6687dd084\\
b75d6f442cc05e7& e265c53617dc284& a19b42c91c823ef& 1db23a9a14842ef\\ 
b8b3709a84250fb& edcfda602e885ae& 47088fee166c2b2& e2b3c49a84250fb\\
1d083bee16cc2a6& ed19dbcfd22fa3a& b7756f142cc05e7& b75d6f47d39fa0c\\ 
edbbda8bd0e7a63& e293c4da84450f7& 47ee8e2215b42c9& a1ef4222e1b5c49\\
0b38178edeedba2& 0b441777b71f69c& 1dee3a2285b50c9& ed39db8cbee97a2\\ 
2f9a5ecbd4c7ae7& f5b3ea982cc05e7& f54deb642cc05e7& d04da167d73fa98\\
a0cf4062866d0b2& fab3f49a179c28c& afb35e9bd0c7a67& 5f30bf9e166c2b2\\ 
85310b9fe937d59& 2fce5e62be6d7b2& f531eb9f4136859& f519ebcf426e832\\
5f30bf9e85350d9& 5f4cbf66873d098& 0b4c1767d73fa98& 0b1817cebd357d9\\ 
d0cfa063d597acd& 2f185fcebd357d9& 2fb25e9abf3d798& 0b641736bcc57e7\\
f5b3ea9b433e818& 05300b9e85350d9& 85e70a3015902cd& a01941ce86cd0a6\\ 
fab3f49a14642f3& 5fb2be9a173c298& a0b3409a173c298& af315f9f419684d\\
a16543341738298& a01941ce85950cd& 5f64bf36173c298& a1654337e867d73\\ 
0d201bbe368c6ae& b0756116c4858ef& b18b62e834206fb& b17563153f7a790\\
e521cbbc9e893ae& 1bf636129e8d3ae& 27204fbe9dd53c5& 27f64e129dd53c5\\ 
e5f7ca13928f22e& 8ddf1a413dd27c5& e48bc8eb94272fb& 4fa29eba34246fb\\
b15d6346c37d810& 8d8b1ae834806ef& e509cbec9d713d1& 8d751b153f7a790\\ 
b08b60ea34246fb& d85db147977f290& 1b5c374794872ef& b0a360bac4258fb\\
1ba236ba9f7d390& b10963ed3e8a7ae& a75d4f449fd9384& 8d751b17cb7f910\\ 
a7a34eb89f79390& 0d5c1b46377c690& 4ff69e1235746d1& a75d4f446f78d90\\
27de4e429dd53c5& 4f5c9f46c4858ef& a7754f146f78d90& e5dfca43928f22e\\ 
192c33a607ac08a& 2b0457f6ce1d9bc& 2bfa560b361e6bc& d4c3a87b36fe6a0\\
199632d207ac08a& 6768cf2e06bc0a8& 811503d46658cb4& 1950335e05a40cb\\ 
d405a9f735e66c3& 81d302599ebb3a8& e6afcca2065c0b4& d487a8f334166fd\\ 
\hline
\end{tabular} \\
\end{center}

\begin{center}
\begin{tabular}{llll}
\multicolumn{4}{c}{Table 2 (continued)} \\  \hline 
d43da98737ee682& 812d03a59ebb3a8& b30567f751e6a43& d441a97ecde59c3\\ 
2368472e8efd1a0& 898712f027a848a& f679ed0e44548f5& f687ecf244548f5\\
2368472e8c151fd& a30547f4ee59db4& c5d38a5b12fe220& 918722f1bebb7a8\\ 
f6ebec2a470c89e& 6f78df0e45448d7& 6feade2a44f48e1& c5698b2f73eee02\\
a3d3465b73eee02& a39746d313ee202& 6f50df5e270c49e& dc69b92f15062df\\ 
09c2127a46bc8a8& a369472f12fe220& 898712f025404d7& 89c3127826b84a8\\
89af12a1bde37c3& c42d89a68d051df& 23d2465a8c151fd& a341477c8f49196\\ 
a3bf46808f49196& a3d3465b13ee202& 918722f3dabfb28& 6fc2de7a24544f5\\
89eb122a40f4861& a369472cec11dfd& f6c3ec7a26bc4a8& dcbfb8828f4d196\\ 
5340a77e270c49e& dc69b92e7ebcfa8& 3bea762be5a7ccb& 5304a7f625e44c3\\
918722f14c129fd& 35fa6a0a24f44e1& 35be6a8245e48c3& 91c322794d029df\\ 
ac3d598644148fd& 9feb3e2b734ee16& cac3947a46fc8a0& 91bf2282b0f5661\\
9f2d3fa773aee0a& c41589d67e5cfb4& 89bf1280d709a9e& 35866af244148fd\\ 
8979130cd7e9a82& acfb580a470c89e& c469892fe547cd7& 9feb3e288da11cb\\
91bf22814de29c3& f92df3a4eeb9da8& 9f153fd7725ee34& c451895fe65fcb4\\ 
f9ebf22b10b6269& dc51b95e7cb4fe9& 9f693f2ced41dd7& f9ebf2288cb11e9\\
9faf3ea08f49196& 9f693f2c8d411d7& 5378a70e44148fd& 9141237cd5e1ac3\\ 
5bc6b67205140dd& 25ba4a8a04e40e3& 8f111fdf524ea36& f193e2d8afb9588\\
814503756ce2de3& 8139038e911525d& 81a302b8979928c& 8fdf1e435196a4d\\ 
81e702316d72dd1& 81a302b96f9ad8c& f14de3675086a6f& f14de367537ea10\\
8fb31e98af79590& 8131039f69d6d45& 81750314979928c& f165e33753dea04\\ 
2bb2569bc4878ef& 81e7023291d5245& 5be6b63205740d1& da5db546073c098\\
d44da967c7df884& 8f091fecad315d9&                &                \\ 
\hline
\end{tabular} \\
\end{center}

\begin{center}
\begin{tabular}{llll}
\multicolumn{4}{c}{Table 3: $H(60)$ from Yang's Theorem 1} \\  \hline 
d9d7bc22eb65e8f& f165ed4722f67bd& e963dd4b21f67dd& 9443270ba1377c5\\
851105aeec3de64& f15bed38c2e9bbe& de17b3a288652ef& f763e14b25f675d\\
82170ba28c25267& de2fb3d1146216f& ef9dd0b4dae98be& 9383288bc237ba5\\
9d2935dd0f22207& c2d78a22e47df6c& 9ad13a2eef3de04& c1178da2eb7de8c\\
831109ae883d2e4& 8f451107cd2fa46& eb73d96a21f47dd& 9be9385d1422167\\
ee9bd2b8c1e9bde& c12f8dd1777ad0c& f65de334a2f17bd& 9bef38517722d07\\ 
9a113bae8f25207& 84d70622883d2e4& 831709a2eb3de84& 9a2f3bd16f3ae04\\
832f09d1773ad04& c6e9825d777ad0c& e85ddf34b9e94de& 88831e8bcd2fa46\\ 
82e90a5d0c22267& c5ef8451647af6c& f69de2b4c2e9bbe& 94452707c237ba5\\
85ef04516c3ae64& f663e34b42eebbe& 9c1137aeeb25e87& 9aef3a510f22207\\ 
8b43190bc22fba6& 8c83168ba12f7c6& d9e9bc5d146216f& 93852887a1377c5\\
e863df4b46eeb3e& 82290bdd6c3ae64& eea5d2c721f67dd& f375e96622f47bd\\ 
f4b3e6ea22f47bd& e865df4725ee75e& f75de134daf18bd& c2e98a5d046236f\\
ecb5d6e621f47dd& dad1ba2ee77df0c& 832909dd143a164& de11b3aeeb65e8f\\ 
c1298ddd147a16c& f775e1664ef4a3d& e875df662dec65e& efa3d0cb25ee75e\\
f0a5eec725f675d& ebb3d8ea41ecbde& 9dd734228f25207& 9d1735a2ef3de04\\ 
c6ef8251147a16c& d9d1bc2e88652ef& 9783208bae37625& f0a3eecb46f6b3d\\
c2178ba2846536f& 8f851087ad37645& daefba51076230f& f19becb8a2f17bd\\ 
90452f07ae37625& c51185aee47df6c& f75be138b9f14dd& c6d78222887d2ec\\
c6d1822eeb7de8c& e9a3dccb41eebde& 84ef0651143a164& f09beeb8daf18bd\\ 
c5d1842e846536f& dde9b45d677af0c& 9c1737a288252e7& e85bdf38dae98be\\
dd29b5dd076230f& e95ddd34c1e9bde& efb3d0ea2dec65e& d9efbc517762d0f\\ 
c2298bdd647af6c& e873df6a4eeca3e& dd17b5a2e77df0c& 9bd73822eb25e87\\
f0b3eeea4ef4a3d& 84e9065d773ad04& f765e14746f6b3d& f3b5e8e642ecbbe\\ 
da2fbbd1677af0c& ec75d76641ecbde& efb5d0e64eeca3e& f773e16a2df465d\\
ddd7b422876530f& f6a3e2cb22f67bd& 9c2937dd7722d07& 8b451907a12f7c6\\ 
ee65d34741eebde& ef9bd0b8b9e94de& f0b5eee62df465d& 9de9345d6f3ae04\\
c1118dae887d2ec& 9743210bce2fa26& c52f85d1046236f& 82d70a22ec3de64\\ 
852f05d10c22267& 85d1042e8c25267& f09deeb4b9f14dd& 9c2f37d11422167\\
90852e87ce2fa26& f473e76a42ecbbe& 8c851687c22fba6& f1a5ecc742eebbe\\ 
88431f0bad37645& de29b3dd7762d0f& efa5d0c746eeb3e& ee5bd338a1f17dd\\
e99ddcb4a1f17dd& 9bd1382e88252e7& 84d1062eeb3de84& da11bbae876530f\\ 
82ed3d3c3350853& 4a32ac809e19dfa& ba034ce32f16bdb& be1344c3ad1fb9a\\
7afccd1f309682b& 6680f5e4be999ea& aa416c676e1e3fa& 72dedd5b61d6203\\ 
9a8f0df8b159812& 668ef5fb309682b& a67f7418a119a1a& 4212bcc02cd0ba3\\
ba0d4cfca119a1a& b63d549c91d9c02& aa4d6c7c0398e4a& 564c947ff39f04a\\ 
86f335033f569d3& 8edf255891d1c03& 7af2cd00be999ea& 92ad1dbc91d1c03\\
a26f7c382310a5b& 92a31da36e163fb& 86fd351cb159812& 92af1db80390e4b\\ 
\hline
\end{tabular} \\
\end{center}

\begin{center}
\begin{tabular}{llll}
\multicolumn{4}{c}{Table 3 (continued)} \\  \hline 
be1d44dc2310a5b& 9e9f05d83350853& 9a810de73f569d3& 564094649e19dfa\\
9e9105c7bd5f992& 4a30ac840c58fb2& a2617c27ad1fb9a& 421cbcdfa2dfa62\\ 
5a7e8c1b20d6a23& 4602b4e0aed9be2& b6315487fc5f1b2& 564294600c58fb2\\
92a11da7fc571b3& 8ed125476e163fb& b63354836e1e3fa& 6290fdc43c909ab\\ 
5a708c04aed9be2& 6eace5bf61d6203& 72dcdd5ff39704b& 4a3eac9bf39f04a\\
6ea2e5a09e11dfb& 460cb4ff20d6a23& 82e33d23bd5f992& 72d0dd449e11dfb\\ 
629efddbb29f86a& b63f54980398e4a& aa4f6c7891d9c02& 8edd255c0390e4b\\
4a3cac9f61de202& 72d2dd400c50fb3& 5e6084242cd0ba3& 7eecc53fb29f86a\\ 
564e947b61de202& 6eaee5bbf39704b& a67174072f16bdb& 7ee2c5203c909ab\\
6ea0e5a40c50fb3& 5e6e843ba2dfa62& aa436c63fc5f1b2& 8ed32543fc571b3\\
\hline
\end{tabular} \\
\end{center}

\begin{center}
\begin{tabular}{llll}
\multicolumn{4}{c}{Table 4: $H(60)$ from Yang's Theorem 2} \\  \hline 
317410f7d9ee45a& 10bcd747b2f4831& 63e084527485d94& 36b41f77c6ee7ba\\
29b2217bbaf6839& 7b26b5de179d1f7& 73a420fbb2f4452& 63217a2e68821eb\\ 
73a411763e979bd& 6ba210f65d8f9bd& c94a2efbb9f1bc5& 08ba52e7d9ee839\\
087ad7465d8f9bd& a5d0f7ce0b9ad88& 17bc6ceba5f67ba& 64208bd26b85e74\\ 
10bce6ca3e975de& 107cd6c65d8f5de& 4228c7c27485d94& 5deec9c26b85217\\
7be74ba20b9ad88& 08bad6c7d1ec831& 4528f84e089de74& 7c274422149ae68\\ 
642145a27782e68& 6ba2217bb2f4831& 087a62ebd9ee45a& 42e8f7ce179dd94\\
7be685d2179dd94& a2d0c8427782e68& 6c621f76428fa5d& 7364217a3e975de\\ 
d14c2f7bb9f17a6& bd16c64268821eb& 10bc5367baf6839& 177d9297b9f1bc5\\
107ce74bb2f4452& 64e175ae778220b& a210f84e149ae68& 087ae6cbd1ec452\\ 
087bac9bc5e9446& 17bcd946428f63e& 5a2ec64274851f7& 31b420fbbaf645a\\
0fbad8c6428fa5d& 0fba6d6bc6ee7ba& 7ce7742e149a20b& 45e8c8426b85e74\\ 
d68c20fba6f1446& 29721177d9ee839& 7ce6ba5e089d217& 0fbba31bdae97a6\\
7c268a52089de74& 107dad1ba6f1446& d64c10f7c5e9446& bad6c9c2778220b\\ 
7b277bae0b9a1eb& 6320b45e74851f7& 64e0bbde6b85217& 08bae74a3e979bd\\
6b621177d1ec831& 0f7b9317dae9bc5& 17bda29bb9f17a6& d18c1f77dae97a6\\ 
107c636bbaf645a& ce4a1177c5e9825& ce8a217ba6f1825& a510c7c26882d88\\
08bb9c97c5e9825& 5aeef64e179d1f7& 6b6220fa5d8f5de& bdd6f64e0b9a1eb\\ 
6ca22f7a428f63e& 5d2ef9ce089d217& 36742f7ba5f67ba& 10bd9d17a6f1825\\
736410f7d1ec452& ba16f9ce149a20b& 63e14a226882d88& c98a1ef7dae9bc5\\ 
1c835c9025753f3& df0f2575e66ea10& 49d7f6383532e04& db9f64717c26a90\\
0073647115a0b96& 713586fd97ae428& aec98fdf2d45ec3& 6dc5bf1ca37bd4d\\ 
ff8d657058b5c4b& fb1d247450b4e82& 923bf63850b4e82& 55a7cfd905e7661\\
b2b9ff1858b5c4b& 96abff19e66f859& 381115b4b73c1ba& 8a5b8ede01d6898\\ 
ff0cd3a83133d4d& e7ecab4d05e7661& fb9d6d517c278d9& aad9cfdbb70de43\\
b23800e515a0b96& aad87800b33d2f3& 8e4a302425741ba& c37f1c972544c0a\\ 
8ada7020b73c1ba& e37d15b7b70cc0a& b2b9fe1850b5ccb& 3c8155b0b73d3f3\\
75a5c7f993af761& aec8390421752f3& fb1d257458b4e02& c7ef5d932d45ec3\\ 
381114b4b33c0ba& e3fcea4997af661& c7ef5c932545e43& b6a9f7397c278d9\\
8a5a7020b33c0ba& e7ed55b3b70de43& b62841e115a19df& e3fceb4993af761\\ 
fb9d6c517427859& 962bb63ccafce02& c3ff54b2939e898& 24612d5583e98df\\
aad87900b73d3f3& 92bbbe1de66ea10& ae49c7fa9b9fa51& c76ea26d01e6528\\ 
db9e9a8ca37af04& 962bb73cc2fce82& ff8cd3a83533c4d& 69d5fe183133d4d\\
df8f6d50c2fdccb& 92ba08c587e8b96& ae498fdf2545e43& 18131c9425741ba\\ 
c76ea36d05e6428& 923a08c583e8a96& 69d5ff183533c4d& 96aa49c187e99df\\
713587fd93ae528& 8ecbceda01d7ad1& c37ee26997ae428& 4d47b73ca37af04\\ 
ff8d2d55e66f859& ff0d657050b5ccb& df8f6c50cafdc4b& 24e12d5587e99df\\
df0f2475ee6ea90& 4d47b63ca77ae04& 51378edd01e6528& 96aa48c183e98df\\ 
b629bf1ccafdc4b& 0463257587e8b96& db1f2c54cafce02& 20f16c5115a19df\\
\hline
\end{tabular} \\
\end{center}

\begin{center}
\begin{tabular}{llll}
\multicolumn{4}{c}{Table 4 (continued)} \\  \hline 
8e4b86ff2544c0a& 3c0155b0b33d2f3& b239b63d7426a10& aa5986fe9b9e818\\ 
923bbe1dee6ea90& 8e4b87ff2d44c8a& db1f64717426a10& b6a9bf1cc2fdccb\\
df8edb88a37bd4d& 7525c7f997af661& 51b78edd05e6428& 96abfe19ee6f8d9\\ 
e3fd5d92939fad1& e76d1d9601d7ad1& db1e9a8ca77ae04& e76d55b3bf0dec3\\
6d45bf1ca77bc4d& b629f7397427859& b6a841e111a18df& e3fd5c929b9fa51\\ 
fb1c92ac3132f04& e76cab4d01e7761& e37d14b7bf0cc8a& ff0d2d55ee6f8d9\\
8ecbcfda09d7a51& aec9c7fa939fad1& c3fee26993ae528& 8e4a312421740ba\\ 
4957f6383132f04& 1c835d9021752f3& 0463247583e8a96& 20f16d5111a18df\\
b23801e511a0a96& aad9cedbbf0dec3& e7ed1d9609d7a51& ae48390425753f3\\ 
8adbc6fbb70cc0a& 8a5bc6fbbf0cc8a& aa5987fe939e898& 92bbf63858b4e02\\
8adb8ede09d6818& 55a7ced901e7761& 18931c9421740ba& 00f3647111a0a96\\ 
c3ff1c972d44c8a& fb1c93ac3532e04& c37f54b29b9e818& c76f15b609d6818\\
df8eda88a77bc4d& c76f14b601d6898& b239b73d7c26a90& db9f2c54c2fce82\\
\hline
\end{tabular} \\
\end{center}

\begin{center}
\begin{tabular}{llll}
\multicolumn{4}{c}{Table 5: $H(60)$ from Yang's Theorem 3} \\  \hline 
3c987fcd5b0aade& 94772e1098f92a0& 1b743015c943897& 8b6b102b893f098\\
fc65fe3538b26a9& 939721d37886eaf& 880b16e8e909c9e& 94692e2ff8cfea6\\ 
ff1bf8c938b26a9& 1c983fcd4942897& a81956ccfb09ede& 8b751014e909c9e\\
b39b61c80ab00e9& d8f7b711aafb4e0& 3b747015db0bade& e7f9c90d5b3aad8\\ 
239641d1b8836af& b38561f76a86cef& c01586d5c973891& c4878ff14972891\\
070408f52aca4e6& 938921ec18b02a9& 039601d1aacb4e6& 881516d7893f098\\ 
b47b6e0beacfce6& ab675030fb09ede& df1bb8c92afa4e0& c36b8029c973891\\
e015c6d5db3bad8& dc65be352afa4e0& aff9590f1b762d1& 90e9272f7886eaf\\ 
a80756f39b3f2d8& b4656e348af90e0& 8ff519146940c97& ac875ff31b762d1\\
c7f9890d4972891& 047a0e092aca4e6& 3fe679315b0aade& 971728d3f8cfea6\\ 
1fe639314942897& 970928ec98f92a0& 180a36e9c943897& afe759307b40ed7\\
270448f538826af& f8f7f711b8b36a9& 247a4e0938826af& 8feb192b0976091\\ 
fb89f1edb8b36a9& b0e567340ab00e9& 8c951fd70976091& ab79500f9b3f2d8\\
db89b1edaafb4e0& 90f7271018b02a9& 380a76e9db0bade& e36bc029db3bad8\\ 
8c8b1fe86940c97& b71b68c88af90e0& b0fb670b6a86cef& 20e8472db8836af\\
e487cff15b3aad8& 00e8072daacb4e6& b70568f7eacfce6& ac995fcc7b40ed7\\
\hline
\end{tabular} \\
\end{center}

\begin{center}
\begin{tabular}{llll}
\multicolumn{4}{c}{Table 6: $H(60)$ from Yang's Theorem 4} \\  \hline 
dfa308e8858c365& 7c74ef6e25dc37a& 2c604f44857176f& d461ba477071490\\
86f51f6dc5dc225& 27a2fdeb658c23a& d6e11a6e30215d0& 27a258c0857176f\\ 
8fb70de9c5dc225& 25225dc090715cf& 753658c1d02148f& 7c744a47708c09a\\
8af70769cfdc365& 22e252408f7162f& 77b6f8ea25dc37a& 2ee04a4490715cf\\ 
2ee0ef6dc52162f& 83b515edcfdc365& 8475bf44858c365& 7536fdea30dc1da\\
7ef4ea6e30dc1da& 2c60ea6dd02148f& dd23a8c37071490& 77b65dc1c52162f\\ 
8d37adc230215d0& 7ef44f47658c23a& 292045c48f7162f& 2522f8eb708c09a\\
72f65741cf2176f& 2ba0e5edcf2176f& dae302688f8c225& 8135b5c48f8c225\\ 
b96fabb29e41fbc& 74dcb0de947cefe& 21881a743ed1bab& ec3b011834ecae9\\
35d832d41e51fbb& 586ee9ba947cefe& 504ef9fac53c4d6& 608c987eb4fcaee\\ 
375837d41a51f3b& 68ac883ee5bc0c6& e41b115865ac0c1& 23081f743ad1b2b\\
b04fb9f2c5014d4& 4c3ec11ab4fcaee& 0fba46103ad1b2b& 0d3a43103ed1bab\\ 
ddf9629ecf01594& 1bea6eb01a51f3b& f86b29b8146cef9& d5d972de9e41fbc\\
dcf9609ec5014d4& 441ed15ae5bc0c6& 808dd87434ecae9& 7cfca09ec53c4d6\\ 
196a6bb01e51fbb& 88adc83465ac0c1& f04b39f8452c4d1& 9cfde094452c4d1\\
94ddf0d4146cef9& b14fbbf2cf01594& 114a7bf04f11593& 3df822944f11593\\
\hline
\end{tabular} \\
\end{center}

\end{document}